\documentclass[11pt]{article}

\usepackage{amsmath, amssymb, amsthm, algorithm, algpseudocode, url, seqsplit, array, booktabs, tabularx, titling}
\usepackage[hidelinks]{hyperref}
\usepackage[
  group-separator = {\,},
  group-minimum-digits = 4
]{siunitx}
\newcolumntype{Y}{>{\raggedright\arraybackslash}X} % ragged right column for tabularx

\title{Constrained PSLQ Search for Machin-like Identities Achieving Record-Low Lehmer Measures}
\author{Nick Craig-Wood \\
  \href{mailto:nick@craig-wood.com}{nick@craig-wood.com} 
}
\date{2025-08-07 \\ {\small Revised 2026-09-16}}

\newtheorem{theorem}{Theorem}[section]
\newtheorem{lemma}[theorem]{Lemma}

\theoremstyle{definition}

\theoremstyle{remark}
\newtheorem{example}[theorem]{Example}
\newtheorem{question}[theorem]{Question}

\newcommand{\acott}[2]{\arctan\left(\tfrac{#1}{#2}\right)}
\newcommand{\acotx}[1]{\arctan\left(\tfrac{1}{#1}\right)}
\newcommand{\acot}[1]{\acotx{\num{#1}}}
\newcommand{\ACOT}[1]{[\seqsplit{#1}]} % compact \acot for big results
\newcommand{\lehmer}[1]{\tag*{\(\scriptstyle \lambda=#1\)}}
\newcommand{\eref}[1]{(\ref{#1})}
\newcommand{\piF}{\tfrac{\pi}{4}}
\newcommand{\WHO}[1]{\hspace*{0em plus 1fill}\makebox{\textbf{(#1)}}}
\newcommand{\omitdigits}[1]{\dots[\num{#1}\,\text{digits omitted}]\dots}

\begin{document}

\maketitle
\vfill
\par
\clearpage

\begin{abstract}
  Machin-like arctangent relations are classical tools for computing $\pi$, with
  efficiency quantified by the Lehmer measure ($\lambda$). We present a framework
  for discovering low-measure relations by coupling the PSLQ integer-relation
  algorithm with number-theoretic filters derived from the algebraic structure
  of Gaussian integers, making large scale search tractable. Our search yields
  new 5 and 6 term relations with record-low Lehmer measures ($\lambda=1.4572,
  \lambda=1.3291$). We also demonstrate how discovered relations can serve as a
  basis for generating new, longer formulae through algorithmic extensions. This
  combined approach of a constrained PSLQ search and algorithmic extension
  provides a robust method for future explorations.
\end{abstract}

\section{Introduction}

Machin-like relations have been used extensively for calculating digits of $\pi$ both in the pre- and post-computer age. The first use was in 1706 when John Machin discovered this formula and then used it to calculate $\pi$ to a then world record 100 decimal places just with pencil and paper.

\begin{alignat}{3}
  \piF &{}={}& 4 & \acot{5} & - & \acot{239} \label{machin} \\
  \lehmer{1.8511}
\end{alignat}

In the computer age the largest $\pi$ computation done with arctan formulae was in 2002. Yasumasa Kanada and his team calculated $\pi$ to 1.24 Trillion digits \cite{Bailey2003} using these two Machin-like formulae (discovered by F. C. W. St\"ormer in 1896 and K. Takano in 1982 respectively)

\begin{alignat}{3}
  \piF &{}={}& 44 & \acot{57}  & +7  & \acot{239} \label{stormer} \\
       &&- 12 & \acot{682} &+ 24 & \acot{12943} \lehmer{1.5860} \\
  \piF &{}={}& 12 & \acot{49}  & +32 & \acot{57} \label{takano} \\
       && - 5 & \acot{239} & +12 & \acot{110443} \lehmer{1.7799}
\end{alignat}

Computational exploration of Machin-like arctangent identities remains a rich source of new formulae for $\pi$ (see \cite{Arndt2001} \cite{Wetherfield1996}
\cite{Wetherfield1997} ).

Historically Todd's process \cite{Wetherfield1996} has been used to find arctan relations. It is a method for ``reducing'' $\acott{a}{b}$ terms, representing them uniquely as sums of irreducible $\acotx{x}$ terms. This process gradually increases the number of terms in the relations, increasing the number of terms in a very elegant way.

In this work we will employ a more brute force strategy using the PSLQ (Partial Sum of Least Squares) algorithm developed by Ferguson and Bailey in 1992 \cite{PSLQ} \cite{Ferguson1999PSLQ} to search systematically for 5 and 6 term relations with low Lehmer measures.

\begin{enumerate}
  \item We describe how to apply PSLQ efficiently to finding arctan relations.
  \item We report new 5 term and 6 term Machin-like formulae with the best-known Lehmer measures to date.
  \item We explain how to extend the relations found by PSLQ into better or longer relations.
\end{enumerate}

This experimental approach not only uncovers novel $\pi$ formulae but also illustrates a general template for using PSLQ when there are many possible relations.

\subsection{Lehmer measure}

In the general form, Machin-like identities take the form

\begin{equation}
\piF = \sum_{n=1}^Nc_n\arctan\big(\frac{a_n}{b_n}\big)
\end{equation}

and their efficiency is often gauged by Lehmer's measure ($\lambda$) \cite{Lehmer1938} which is defined as

\begin{equation}
\lambda
=\sum_{n=1}^{N}
  \frac{1}{\log_{10}\!\bigl(\tfrac{b_n}{a_n}\bigr)}
\end{equation}

which is proportional to the total number of $\arctan$ terms to compute $\pi$ to a given accuracy. A lower Lehmer measure is better so it is a useful proxy for efficiency.

In this paper we've concentrated on relations where $a_n=1$ as it is very easy to make low Lehmer measure relations with that constraint relaxed.

\subsection{Current best known relations} \label{sec:best}

 Many Machin-like formulae have been discovered over the last 300 years. Gauss in particular was obsessed with finding them. The relations with the smallest known Lehmer measure for a given number of terms are detailed below.

We have seen the best relation for 2 terms already, it is Machin's original formula \eref{machin}.

The best relation for 3 terms is this one by Gauss (published posthumously in 1863)

\begin{alignat}{4}
  \piF &{}={}& 12 & \acot{18} & +8 & \acot{57} & -5 & \acot{239} \\
  \lehmer{1.7866}
\end{alignat}

and the best relation for 4 terms we have already seen as St\"ormer's \eref{stormer}. Until the relatively recently this was the best known relation.

The previous record-holder for a 5 term relation was found by Amrik Singh Nimbran in 2009 and published on Wetherfield and Chien-Lih's Machination website \cite{Machination} (preserved in the machin-like.org database \cite{MachinLike})

\begin{alignat}{3}
  \piF &{}={}& 22 & \acot{28} & +2 & \acot{443} \label{nimbran5} \\
       && -5 & \acot{1112} & -5 & \acot{8637069} \notag \\
       && +5 & \acot{133280317780182} \lehmer{1.6121}
\end{alignat}

which improved on St\"ormer's 5 term relation from 1896

\begin{alignat}{3}
  \piF &{}={}& 88 & \acot{172} & +51 & \acot{239} \label{stormer5} \\
       && +32 & \acot{682} & +44 & \acot{5357} \notag \\
       && +68 & \acot{12943}. \lehmer{1.7320}
\end{alignat}

The best 6 term relation in the journal literature is this one found by Hwang Chien-Lih in 1997 \cite{Chien1997}

\begin{alignat}{3}
  \piF &{}={}& 183& \acot{239}     &+32  & \acot{1023} \label{chien} \\
       &&  -68& \acot{5832}    & +12 & \acot{110443} \notag  \\
       &&  -12& \acot{4841182} &-100 & \acot{6826318} \lehmer{1.5124}
\end{alignat}

however Wetherfield found a better one in 2004 \cite{Machination}. It was published as a 5 term relation with one term of the form $\acott{2}{x}$; expanding that term into two $\acotx{x}$ terms with \eref{twoa} gives the previous best 6 term relation

\begin{alignat}{3}
  \piF &{}={}& 83 & \acot{107}    & +17 & \acot{1710} \label{wetherfield6} \\
       &&  -22 & \acot{105218} & -24 & \acot{2513489} \notag \\
       &&  -22 & \acot{7167807} & +12 & \acot{7939642926390344818}. \lehmer{1.3562}
\end{alignat}

The best known 7 term formulae was discovered in the same way by Wetherfield in 2004 \cite{Machination}

\begin{alignat}{2}
  \piF &{}={}&  83 & \acot{107} \, +17 \, \acot{1710} \label{wetherfield7}\\
       &&  -22 & \acot{103697} \, -24 \, \acot{2513489} \notag\\
       &&  -44 & \acot{18280007883} \, +12 \, \acot{7939642926390344818} \notag\\
       &&  +22 & \acot{3054211727257704725384731479018} \lehmer{1.3408}
\end{alignat}

\section{How to apply PSLQ to finding arctan relations}

PSLQ finds integer relationships between real numbers. If given a set of $y_i \in \mathbb{R}$, the PSLQ algorithm finds $m_i \in \mathbb{Z}$ such that $\sum m_i y_i = 0$.

This means we can set $y_0 = \pi$ and $y_i = \acotx{x_i}$ for some set of $x_i \in \mathbb{Z}$ and if we are lucky it will produce an integer relation with some $\acotx{x_i}$ terms and $\pi$ thus finding an Machin-like relation for us.

PSLQ is guaranteed to find a relation (if one exists), but it may not find the relation we want. It is quite likely to return relations which do not include $\pi$ as there are many arctan relations which sum to zero rather than $\pi$. A simple example is

\[
\acot{3} - \acot{5} - \acot{8} = 0.
\]

It is also quite likely to return relations which have a high Lehmer measure which we aren't interested in, or relations which sum to multiples of $\pi$ which can represent linear combinations of other relations. In order to fix this we will have to do a combinatorial search only providing a subset of the $\acotx{x_i}$ terms and run PSLQ many times. This is most easily done by masking out terms using a bit mask created using Gosper's Hack \cite[p. 22]{HackersDelight}. Unfortunately this is $O(n!)$.

Though the running time of PSLQ is polynomial the order of the polynomial is quite high, something like $O(n^5)$ \cite{PSLQ} multiplies in the worst case.

Given the run time of PSLQ ($O(n^5)$) and the fact we will have to do a combinatorial search as well ($O(n!)$), this means we want to keep the number of $x_i$ to a minimum. This in turn means we need a clever method to select the $x_i$ to ensure we don't have too many and to maximise our chances that there could be solutions within the $x_i$. An exhaustive approach of setting $x_i$ to all the integers less than $N$ just won't work in a reasonable time.

In the next few sections we'll develop a bit of theory to help us reduce the number of $x_i$ we need and methods to filter those before even running PSLQ.

\subsection{Prime Divisors of \texorpdfstring{\(x^2+1\)}{x\^2+1}}

We'll need some results on the prime divisors of \(x^2+1\) so let's get those out of the way.

\begin{lemma}
For any integer $x$, if a prime $p$ divides $x^2+1$, then $p=2$ or $p\equiv1 \pmod 4$. \label{normprimes}
\end{lemma}

\begin{proof}
First suppose $p>2$ and $p\mid x^2+1$.  Then using Euler's criterion
\begin{align*}
x^2 &\equiv -1 \pmod p
\;\Longrightarrow\;
x^{p-1} \equiv (-1)^{\frac{p-1}{2}} \pmod p.
\end{align*}
By Fermat’s little theorem $x^{p-1}\equiv1\pmod p$, so
\[
(-1)^{\frac{p-1}{2}} = 1
\;\Longrightarrow\;
\frac{p-1}{2}\text{ is even}
\;\Longrightarrow\;
p\equiv1\pmod4.
\]
Finally, suppose $p = 2$, if $x$ is odd then $2\mid x^2+1$ completing the proof.
\end{proof}

\begin{lemma}
Let $x\in\mathbb{Z}$. If \(2\mid x^2+1\) then \(4\nmid x^2+1\). \label{normtwos}
\end{lemma}

\begin{proof}
  We consider two cases.  
  If $x$ is even then $x^2+1\equiv1\pmod4$ so $4\nmid x^2+1$.
  If $x$ is odd, write $x=2k+1$.  Then
  \[
  x^2+1 =(2k+1)^2+1 =4k^2+4k+2 \equiv2\pmod4,
  \]
  so again $4\nmid x^2+1$.
\end{proof}

\subsection{Artcan and the Gaussian integers}

Machin-like relations can be related quite easily to complex numbers and as we are only interested in rational arguments to the $\arctan$ function, to the Gaussian integers \(\mathbb{Z}[i]\).

\subsubsection{Machin-like identities via complex numbers}

Assuming some basic results from complex analysis:

\begin{align*}
\arg(x+i)     &=  \arctan(\tfrac1x) \\
\arg(z_1 z_2) &=  \arg(z_1) + \arg(z_2) \\
\arg(z^n)     &=n \arg(z).
\end{align*}

A Machin-like formula

\[
\piF=\sum_k a_k\arctan\frac1{x_k}
\]

is therefore equivalent to

\begin{align*}
\arg(1+i)
&=\sum_k a_k\arg({x_k}+i) \\
&=\sum_k \arg(({x_k}+i)^{a_k}) \\
&=\arg(\prod_k ({x_k}+i)^{a_k}). \\
\end{align*}

If $\arg(z_1) = \arg(z_2)$ then $z_1$ and $z_2$ must be real multiples
of each other. However we know the $x_k$ and $a_k$ are integers so in
this case they must be rational multiples

\[
u(1+i)= \prod_k ({x_k}+i)^{a_k}, \quad u \in \mathbb{Q}.
\]

For example we can use this to prove Machin's formula \(\piF = 4\,\acot{4} - \acot{239}\). Substituting the RHS of Machin's formula into the equation above and simplifying we get

\begin{align*}
  (5+i)^{4} \cdot (239+i)^{-1} &= (476+480i) \cdot (239-i)/57122 \\
  &= (14244 + 114244\,i) /57122 \\
  &= 2(1+i).
\end{align*}

Since this is of the form expected (a rational $u=2$ multiplied by $1+i$) Machin's formula is must be correct.

\subsubsection{Gaussian integers and constraints on the norms \(N_k=x_k^2+1\)}

For a Gaussian integer $z = a + bi$ where $a,b\in\mathbb{Z}$, its norm, denoted as $N(z)$, is defined as

\[
N(z) = a^2 + b^2.
\]

This is always an integer. The norm function has some useful properties such as

\[
N(z_1 z_2) = N(z_1) N(z_2).
\]

If $N(z)$ is a prime then $z$ is a Gaussian prime. An ordinary prime \(p\in\mathbb{Z}\) has one of three behaviours in \(\mathbb{Z}[i]\):

\begin{itemize}
  \item \emph{Ramified} (\(p=2\)): 
    \[
      2=-i(1+i)^2,\quad N(1+i)=2.
    \]
  \item \emph{Split} (\(p\equiv1\pmod4\)): write \(p=a^2+b^2\). Then
    \[
      p=(a+bi)(a-bi),\quad N(a\pm bi)=p.
    \]
  \item \emph{Inert} (\(p\equiv3\pmod4\)): remains prime in \(\mathbb{Z}[i]\), with
    \[
      N(p)=p^2.
    \]
\end{itemize}

Any Gaussian prime factor $z$ of $x_k+i$ must have $N(z)$ dividing $N(x_k+i) = x_k^2+1$. As shown in lemma~\ref{normprimes} the only prime factors $x^2+1$ can have are 2 or \(p\equiv1\pmod4\). This means Gaussian factors of $x+i$ can be \emph{ramified} and \emph{split} but not \emph{inert}. For the \emph{split} primes

$$
p=(a+bi)(a-bi), \quad N(a\pm bi)=p.
$$

If $p\mid x_j^2+1$, then $a+bi\mid(x_j+i)$ or $a-bi\mid(x_j+i)$ (but not both as that would make an integer). Without loss of generality take $a+bi\mid(x_j+i)$. Each $x_k+i$ will have a unique (up to units) factorization into Gaussian primes which may include 2 and \emph{split} primes, but no \emph{inert} primes. Recalling our complex number version of Machin-like relations

\begin{equation}
u(1+i)= \prod_k ({x_k}+i)^{a_k} \label{zmachin}
\end{equation}

we can see the only Gaussian prime dividing the left of \eref{zmachin} is the \emph{ramified} prime $1+i$, so all the \emph{split} primes on the right of \eref{zmachin} must cancel.

But if $a+bi$ divides \emph{only one} $(x_j+i)$ its exponent would be $a_j\ne0$ which would mean this Gaussian prime would appear on the left of \eref{zmachin}. Therefore it must be cancelled out by either $a+bi$ in another term with a negative $a_j$ or $a-bi$ in a term with a positive $a_j$ or a combination of both. Therefore some other $x_{k\ne j}+i$ must be divisible by $a\pm bi$ and hence some other $x_{k\ne j}^2+1$ is divisible by p.

From this we can conclude that if $x_j^2+1$ has a prime factor other than 2, then this prime factor must appear in at least one other term.

This gives us a powerful filter for selecting $x_k$ to apply PSLQ to.

\subsection{Constraints on \texorpdfstring{$x_i$}{x\_i}} \label{conditions}

We now have developed a set of constraints on the prime factors of the norms $x_i^2+1$ in an arctan relation:

\begin{itemize}
\item Any prime factor $p \mid x_i^2+1$ must be either $2$ or a prime $p \equiv 1\pmod{4}$
\item If $x_i^2+1$ is divisible by 2 then it won't be divisible by 4 or any higher power of 2.
\item If $p \mid x_i^2+1$ where $p \equiv 1\pmod{4}$ then for at least one other $x_j,\;j\ne i$, $p \mid x_j^2+1$.
\end{itemize}

\subsection{Selection of the \texorpdfstring{$x_i$}{x\_i}}

We can now create a set of $x_i$ which have a good chance of creating an arctan relation starting with a set of primes P such that $p_i = 1\pmod4$ . Assuming we are looking for $x_i$ less than some limit $N$ we create a set of all possible numbers $C < N$ which obey the above constraints. This is easy to do with a depth first tree traversal algorithm (Algorithm~\ref{psmooth}).

\begin{algorithm}
\caption{Enumerate candidates with factors from $P$ up to $N$}
\label{psmooth}
\begin{algorithmic}[1]
\Procedure{DFS}{$i,\,c$}
  \Comment{Start with \Call{DFS}{0,\,1}}
  \If{$i = \mathrm{length}(P)$}
    \State \Call{CheckCandidate}{$c$}
    \If{$c \times 2 < N$}
      \State \Call{CheckCandidate}{$c \times 2$}
    \EndIf
    \State \Return
  \EndIf
  \State \Call{DFS}{$i+1,\,c$}
  \State $p \gets \texttt{P}[i]$
  \State $tmp \gets c$
  \While{$tmp \times p < N$}
    \State $tmp \gets tmp \times p$
    \State \Call{DFS}{$i+1,\,tmp$}
  \EndWhile
\EndProcedure
\end{algorithmic}
\end{algorithm}

The candidates $c$ will have the form $\prod_k p_k^{e_k}$ or $2\,\prod_k p_k^{e_k}$ where $e_k \in \mathbb{Z}_{\ge0}$. We then attempt to solve for $x^{2}+1=c$ where $c\in C$. If $\sqrt{c-1} \in \mathbb{Z}$ then we have found an $x_i$. We store the $x_i$ and the factorization of $x_i^2+1$. Since solutions are rare it is quicker to factor $c$ again here using the small set of known primes in $P$, rather than keep track of the factorization in the depth first search (Algorithm~\ref{checkCandidate}).

\begin{algorithm}
\caption{Check if candidate $c$ is a solution of $x^2+1$}
\label{checkCandidate}
\begin{algorithmic}[1]
\Procedure{CheckCandidate}{$c$}
  \State $d \gets c - 1$
  \If{$d < 2$}
    \State \Return
  \EndIf
  \State $t \gets \lfloor \sqrt{d} \rfloor$
  \If{$t \times t = d$}
    \State $x \gets t$
    \State $factors \gets \textproc{FactorInteger}(c,\,P)$
    \State \textproc{RecordSolution}($x,\,factors$)
  \EndIf
\EndProcedure
\end{algorithmic}
\end{algorithm}

Once we have found all possible $x_i < N$ then we check to see whether each $p\in P$ appears twice or more in the factorizations of the $x_i^2+1$. If it does not then we can discard this set of $x_i$ as they cannot make an arctan relation using all the primes in $P$.

If the $x_i$ do pass this test then we can run PSLQ on the $\acotx{x_i}$ and $\pi$ to attempt to find an arctan relation. We then repeat this on all subsets of $x_i$ making sure the conditions in \ref{conditions} above still hold for every subset.

\subsection{Selection of \texorpdfstring{$P$}{P}}

We now have an effective method of finding $x_i$ given a set of primes $P$. The remaining question remaining is how to choose the set of primes $P$ such that it produces good arctan relations. A number of possible strategies were tried:

\begin{itemize}
\item Use a straight forward combinatorial approach choosing $i$ primes $p \equiv 1\pmod{4}$ less than a limit and using these as $P$. This works well when $i$ was small but as $i$ got much above 6 the runtimes became too long. This method was successfully used to to rediscover many the existing arctan relations including \eref{chien}.
\item Use the unique primes in the factorization of $x^2+1$ for $x<10^8$ (excluding $2$) as $P$. This worked surprisingly well and it quickly found new relations, including the best 5 term relation and some very good 6 term relations.
\item Evolve a set of primes. It was noticed that the set of primes used by existing low Lehmer measure relations are often based on a previous not so high performing relation with one extra prime. A strategy was developed which kept a few hundred of the highest performing prime sets and tried to add one more prime to each set, measured the results and trimmed the set back to the few hundred best performers. This turned out to be very successful and is how the best 6 term relation was discovered.
\end{itemize}

\subsection{Computational tricks}

In order to speed up the computation various tricks were used.

\begin{itemize}

\item The combinatorial search through the PSLQ was only allowed to run for a threshold number of iterations without finding a new result in order to abandon prime groups which had been ``worked out''. The combinatorial search was also abandoned if the PSLQ failed for too many iterations in a row - this generally indicated that there were no more relations to be found.

\item When the combinatorial search found a relation (a true relation or a zero relation) the bit mask for this was stored and any combinations with this set of inputs was skipped. This heuristic was found to skip through the search much quicker as generally PSLQ seems to favour results it had already found.

\item When doing the combinatorial search, when items were pruned, the PSLQ was skipped if the conditions in \ref{conditions} did not hold. This pruned the search massively.

\item When searching for longer relations, these are much more likely to start with a higher initial $\acotx{x}$ term. Since there are many solutions when x is small, discarding $x$ if it was less than some $\alpha$ was very effective at pruning the number of $x_i$ in the search.

\end{itemize}

\section{Results}

We found $\num{708024}$ Machin-like arctan relations with a Lehmer measure better than $1.86$ ranging from 2 to 15 terms. We also found many millions more relations with a Lehmer measure worse than $1.86$ but we discarded those as our focus was finding relations with low Lehmer measure.

The most interesting relations were these two being respectively the 5 and 6 term relations with the lowest Lehmer measure known (to the best of our knowledge) for arctan relations of the form $\acotx{x}$.

\begin{alignat}{3}
  \piF &{}={}& 29 & \acot{68}       & +42 & \acot{117}      \label{ncw5} \\
       && -15 & \acot{2675143}  & -5  & \acot{2976163}  \notag\\
       &&  -5 & \acot{302342643} \lehmer{1.4572} \\
  \piF &{}={}& 83 & \acot{107}      &  +17 & \acot{1710}    \label{ncw6}\\
       && -44 & \acot{225443}   & -68 & \acot{2513489}  \notag\\
       && +22 & \acot{42483057} & +34 & \acot{7939642926390344818} \lehmer{1.3291}
\end{alignat}

\eref{ncw6} is now the relation with the lowest known Lehmer measure superseding \eref{wetherfield7} (excluding relations with very large numbers of digits such as \eref{alferov}).

A table of the best relations found for each number of terms can be seen in table~\eref{table:terms1} and table~\eref{table:terms2}.

We searched the literature (for example \cite{Arndt2001} \cite{Chien1997} \cite{Machination}) and the Internet extensively looking for published Machin-like relations. The most complete source is the machin-like.org database \cite{MachinLike}, which incorporates the relations from the Machination website \cite{Machination}. Many of those relations contain half-integer terms $\acott{2}{x}$; when these are expanded into two $\acotx{x}$ terms using \eref{twoa} the only published relation better than those in our tables is Wetherfield's 9 term relation from 2009, which we show in place of ours in table~\eref{table:terms1}.

\subsection{Example relations}

Here are some examples of $P$, the $x_i$ found and the arctan relations PSLQ with combinatorial search finds.

\begin{example}

Prime group $P=\{5\}$ gives $x_i = \{2,3,7\}$. This gives the classical arctan formulae which are attributed to Hutton (1776), Herman (1706) and Euler (1738) though it is now thought that Machin probably discovered them all first.

\begin{alignat*}{3}
\piF &{}={}& 2 & \acot{3} &+& \acot{7} \\
\piF &{}={}& 2 & \acot{2} &-& \acot{7} \\
\piF &{}={}&   & \acot{2} &+& \acot{3}.
\end{alignat*}

\end{example}

\begin{example}

Prime group $P=\{13\}$ has $x_i = \{5, 239\}$ which gives Machin's formula \eref{machin}.

\begin{alignat}{3}
\piF &{}={}& 4 \acot{5} - \acot{239}. \\
\lehmer{1.8511}
\end{alignat}

\end{example}

\begin{example}

Prime group $P=\{5, 13\}$ has $x_i = \{2,3,5,7,8,18,57,239\}$
which gives 42 relations, the best of which was discovered by Gauss and published in 1863

\begin{alignat}{3}
\piF &{}={}& +12 \acot{18} +8 \acot{57} -5 \acot{239}. \\
\lehmer{1.7866}
\end{alignat}

\end{example}

\begin{example}

Using
\begin{alignat*}{2}
P   &{}={}& \{&5,13,61\} \\
x_i &{}={}& \{&2,3,5,7,8,11,18,57,239,682,12943\}
\end{alignat*}

gives 111 relations, the best of which is by St\"ormer's 4 term \eref{stormer}

\begin{alignat}{3}
  \piF &{}={}& 44 &\acot{57}  & +7& \acot{239} \notag\\
       &     &-12 &\acot{682} &+24& \acot{12943}. \lehmer{1.5860}
\end{alignat}

\end{example}

\begin{example}

Using
\begin{alignat*}{2}
P   &{}={}& \{&5, 13, 229, 457, 1201\} \\
x_i &{}={}& \{&2,3,5,7,8,18,49,57,107,109,122,239, \\
    & &   &1023,5832,110443,4841182,6826318\}
\end{alignat*}

we find 1216 relations, with 3 relations better than St\"ormer's formula the first of which was found by Hwang Chien-Lih in 1997 \cite{Chien1997} and was, until Wetherfield's 2004 relation \eref{wetherfield6}, the 6 term relation with the lowest known Lehmer measure. You can see that these relations have many identical arctan terms and it is common when using this method to find lots of very similar arctan relations. Note that the last two relations don't sum to $\piF$.

\begin{alignat}{4}
  \piF &{}={}& 183& \acot{239}    &+32& \acot{1023}    &-68& \acot{5832} \\
       & & +12& \acot{110443} & -12& \acot{4841182} &-100& \acot{6826318} \lehmer{1.5124} \\
  \pi &{}={}& 366& \acot{122}    &+128& \acot{1023}    &+94& \acot{5832} \\
      & & +48& \acot{110443} & -48& \acot{4841182} &-34& \acot{6826318} \lehmer{1.5713} \\
  \pi &{}={}& 366& \acot{109}    &-238& \acot{1023}    &+94& \acot{5832} \\
      & & +48& \acot{110443} & -48& \acot{4841182} &-34& \acot{6826318} \lehmer{1.5828}
\end{alignat}

\end{example}

\begin{example}

%% ./find_arctan_formulae -skip-relations -group "5,13,37,24113,76369"

Using

\begin{alignat*}{2}
P   &{}={}& \{&5,13,37,24113,76369\} \\
x_i &{}={}& \{& 2,3,5,6,7,8,18,31,43,57,68,117,239,2228,\\
    & &   &1408818,2675143,2976163,302342643\}
\end{alignat*}

produces 1527 relations the best of which are 3 clearly related 5 term relations. The first is \eref{ncw5} which is the best 5 term relation found.

\begin{alignat}{3}
  \piF &{}={}& 29& \acot{68}      &+42& \acot{117} \\
       &&  -15& \acot{2675143} & -5& \acot{2976163} \notag\\
       &&  -5& \acot{302342643} \lehmer{1.4572} \\
  \piF &{}={}& 29& \acot{68}      &+42& \acot{117}  \\
       &&  -15& \acot{1408818} & +10& \acot{2976163} \notag\\
       &&  -5& \acot{302342643} \lehmer{1.4642} \\
  \piF &{}={}& 29& \acot{68}      &+42& \acot{117} \\
       &&  -5& \acot{1408818} & -10& \acot{2675143} \notag\\
       &&  -5& \acot{302342643} \lehmer{1.4654}
\end{alignat}

\end{example}

\subsection{Extending the search with \texorpdfstring{$\acott{2}{a}$}{arctan(2/a)} terms}

It was noticed that it was possible to produce relations of much lower Lehmer measure if the constraint of using $\acotx{a}$ terms was relaxed. If terms of the form $\acott{2}{a}$ were allowed then we could produce relations like this with very low Lehmer measure (the first of which was discovered first by Wetherfield in 2004).

\begin{alignat}{3}
  \piF  &{}={}& 83 & \acot{107}    & +17 & \acot{1710} \notag\\
        && -22 & \acot{103697} & -12 & \acott{2}{\num{2513489}} \notag\\
        && -22 & \acott{2}{\num{18280007883}} \lehmer{1.2657} \\
  \piF  &{}={}& 83 & \acot{107}    & +17 & \acot{1710} \notag\\
        && -44 & \acot{225443} & -34 & \acott{2}{\num{2513489}} \notag\\
        && +22 & \acot{42483057} \label{ncw6precursor} \lehmer{1.2839}
\end{alignat}

%% ./find_arctan_formulae -skip-relations -max-numerator 2 -group "5,113,229,177553"

The second was discovered with this $P$ and $x_i$

\begin{alignat*}{2}
P   &{}={}& \{&5,113,229,177553\} \\
x_i &{}={}& \{&2,3,7,11/2,15,15/2,107,1710,,\\
    & &   &225443,2513489/2,42483057\}
\end{alignat*}

We can then expand out the $\acott{2}{a}$ terms using this formula which was modified from Wetherfield's 1996 work on Todd's process \cite{Wetherfield1996}. It is easy to prove with repeated applications of \eref{addarctans}.

\begin{equation}
\acott{2}{a} = 2\,\acotx{a} + \acott{-2}{a^3 + 3a} \label{twoa}
\end{equation}

Note that $a^3 + 3a$ is always divisible by 2 if $a \in \mathbb{Z}$ which means that this will always produce two $\acott{1}{x}$ terms.

These two relations above, when the $\acott{2}{a}$ terms are expanded out with \eref{twoa} make the best 7 term formula \eref{wetherfield7} and make the best 6 term formula \eref{ncw6} respectively.

No extra improvement was found using $\acott{3}{a}$ terms. Adding more terms means that there are more possible relations to choose from, which in turn puts more pressure on the PSLQ and combinatorial search.

\subsection{Extending the search with the floor-function iteration} \label{sec:iteration}

In 2022 Abrarov, Siddiqui, Jagpal and Quine \cite{Abrarov2022} introduced an iteration which turns the rational remainder of a partial Machin-like relation into a sequence of $\acotx{x}$ terms with integer $x$, using the floor function to pick each $x$, and observed that the integers grow very rapidly so the Lehmer measure stays small (see also \cite{Abrarov2023}). In 2023 Oleg S. Alferov \cite{Alferov2023} developed this into a process which starts from any term $m \acotx{q_0}$, chooses the nearest integer at each step, always terminates, and creates a Machin-like formula of an arbitrarily small Lehmer measure.

The process is quite simple. Using a starting point term $m \acotx{q_0}$, we will add or subtract $\acotx{x}$ terms using this trigonometric identity reducing the difference between this and $\piF$ each time, eventually to zero

\begin{equation}
\arctan(x) + \arctan(y) = \arctan\left(\frac{x + y}{1 - xy}\right).  \label{addarctans}
\end{equation}

We can use this identity to multiply an arctan term by an integer by repeated addition, for example

\begin{equation}
2 \arctan(x) = \arctan\left(\frac{2 x}{1 - x^2}\right). \label{mularctans}
\end{equation}

This can easily be continued to any multiple of $\arctan(x)$.

In Alferov's paper he starts with a term $m \acotx{q_0}$ where m is chosen to make $m \acotx{q_0}$ as close to $\arctan(1) = \piF$ as possible. The continued fractions of $\piF$ are ideal for this.

We first calculate the remainder such that $\arctan{r} = \arctan(1) - m \acotx{a}$. We do this by repeated use of \eref{addarctans} and \eref{mularctans} to turn the $m \acotx{q_0}$ into $\acott{a}{b}$ then using \eref{addarctans} again to add $-\acott{a}{b}$ to $\arctan(1)$. This will get us a rational $r$ value. This is the starting point for the iteration and it should be as small as possible.

We then reduce the remainder $\arctan r = \acott{a_{n}}{b_{n}}$ by subtracting $\arctan\left(1 / \left\lfloor \frac{b_n}{a_n} \right\rfloor \right)$ again using \eref{addarctans}. The remainder is guaranteed to become zero in a finite number of terms provided $|b_0| > |a_0| > 0$

\begin{proof}
  \begin{align*}
    \acott{a_{n+1}}{b_{n+1}}
    &= \acott{a_n}{b_n} - \acotx{\left\lfloor \frac{b_n}{a_n} \right\rfloor} \\
    &= \acott{a_n \left\lfloor{\frac{b_n}{a_n}}\right\rfloor - b_n}{a_n + b_n \left\lfloor{\frac{b_n}{a_n}}\right\rfloor} && \text{using \eqref{addarctans}} \\
    &= \acott{a_n \left({\frac{b_n}{a_n}} - \frac{b_n \bmod a_n}{a_n}\right) - b_n}{a_n + b_n \left\lfloor{\frac{b_n}{a_n}}\right\rfloor} \\
    &= \acott{-\left(b_n \bmod a_n\right)}{a_n + b_n \left\lfloor{\frac{b_n}{a_n}}\right\rfloor} \\
  \end{align*}
  In the worst case (with no cancellations with $b_{n+1}$), $a_{n+1}$ is $-b_n \bmod a_n$ so $|a_{n+1}| \le |b_n \bmod a_n| < |a_n|$. Since $a_n$ is an integer and $|a_{n+1}| < |a_n|$ the sequence of $|a_n|$ values strictly decreases and must come to an end with $a_k = 0$.
\end{proof}

If we assume no cancellations between the $a_n$ and $b_n$ then on average we will halve $a_n$ each iteration which gives a rough estimate of the number of terms as $\log_2{a_0}$. In practice there will likely be fewer terms than this as $\gcd(a_n, b_n)$ will be greater than one and sometimes much greater.

Alferov proved for this process

\[
\lambda < \frac{3}{\log_{10} q_0}
\]

which means that by selecting a large enough $q_0$ we can make a finite sequence with as low a Lehmer measure as desired.

The above is a slight simplification of the method - in practice we may choose $\arctan\left(1 / \left\lfloor \frac{b_n}{a_n} \right\rfloor \right)$ or $\arctan\left(1 / \left( \left\lfloor \frac{b_n}{a_n} \right\rfloor + 1 \right) \right)$ whichever gives the lowest absolute remainder. We also found that selecting for the smallest $a_{n+1}$ after cancellation instead of choosing the lowest absolute remainder was sometimes advantageous in producing shorter sequences.

There is a practical problem with this method though. It can be seen above that the $b_{n+1} \approx b_n^2/a_n$. Since $a_n$ is small and reducing, this means that the $b_n$ roughly double in number of digits each iteration. This in turn means that for quite small $q_0$ we get $b_n$ with millions of digits. For example starting with $22 \acot{28}$ gives us the 23 term sequence below with a Lehmer measure of $1.0919$ which is very low, but unfortunately $b_{22}$ has $\num{11512147}$ digits!

\allowdisplaybreaks             % allow pages to break equations

\begin{alignat}{3}
  \piF &{}={}& 22 & \acot{28} \label{alferov}\\
&& + & \acot{56547} \notag\\
&& + & \acot{20747394343} \notag\\
&& + & \acot{1112172624652580034840} \notag\\
&& - & \acot{16659543628852678157467292276729792021493732} \notag\\
&& + & \acotx{19351\omitdigits{78}88894} \notag\\
&& + & \acotx{14718\omitdigits{166}74124} \notag\\
&& - & \acotx{53767\omitdigits{342}09023} \notag\\
&& - & \acotx{12284\omitdigits{693}30646} \notag\\
&& + & \acotx{38904\omitdigits{1395}93265} \notag\\
&& + & \acotx{86261\omitdigits{2801}59064} \notag\\
&& - & \acotx{29087\omitdigits{5611}59141} \notag\\
&& + & \acotx{17858\omitdigits{11232}65666} \notag\\
&& - & \acotx{15170\omitdigits{22475}95919} \notag\\
&& - & \acotx{56421\omitdigits{44959}94705} \notag\\
&& - & \acotx{67328\omitdigits{89928}05021} \notag\\
&& - & \acotx{12422\omitdigits{179867}42982} \notag\\
&& + & \acotx{48665\omitdigits{359743}41495} \notag\\
&& + & \acotx{20525\omitdigits{719499}86918} \notag\\
&& - & \acotx{37353\omitdigits{1439008}30567} \notag\\
&& - & \acotx{83311\omitdigits{2878027}19649} \notag\\
&& - & \acotx{16923\omitdigits{5756063}78545} \notag\\
&& - & \acotx{17656\omitdigits{11512137}79193} \lehmer{1.0919}
\end{alignat}

\interdisplaylinepenalty=10000  % forbid pages to break equations

We can, however, use the iteration to create longer versions of existing relations. (Abrarov et al.\ \cite{Abrarov2022} used it in a similar spirit, to convert the half-integer terms of Wetherfield's 2004 relation \eref{wetherfield6} into integer terms.) If we truncate the existing relation we will likely have a very small residue where the process works best. We can use this as the initial term and can complete the relation using the iteration. This may or may not improve the Lehmer measure but it does give interesting relations with more terms. If we see that any remaining terms had a common multiplier then we can use this to reduce the initial term further. This is how the best relations with 8 terms and above were found (arbitrarily discarding sequences with $b_n > 2^{1000}$).

For example the 8 term sequence with the smallest Lehmer measure found

\begin{alignat}{2}
  \piF &{}={}& 83 & \acot{107} \label{ncw8}\\
       && +17 & \acot{1799} \notag\\
       && +29 & \acot{110443} \notag\\
       && +22 & \acot{4841182} \notag\\
       && +17 & \acot{5675477} \notag\\
       && -17 & \acot{230921399215798} \notag\\
       && -17 & \acot{180288246462881610746480172447} \notag\\
       && +17 & \acot{68258088806587561619922885741255927889500440723402481904848} \notag\\
  \lehmer{1.4167}
\end{alignat}

was found by truncating this formula (which was found with PSLQ) after the 4th term and using the iteration to complete it, improving the Lehmer measure slightly and increasing the number of terms by 1.

\begin{alignat}{2}
  \piF &{}={}& 83 & \acot{107} \\
       && +17 & \acot{1799} \notag\\
       && +29 & \acot{110443} \notag\\
       && +22 & \acot{4841182} \notag\\
       && +17 & \acot{5724143} \notag\\
       && +17 & \acot{667557142} \notag\\
       && +17 & \acot{1364321153627926317} \lehmer{1.4643}
\end{alignat}

After the PSLQ process had found as many relations as we thought it was going to we tried this process on all the relations, truncating at all possible places.

\subsection{Computation}

The computations were done with $x < 2^{64}$ which keeps $x^2+1 < 2^{128}$ which is convenient for a 64 bit architecture. The maximum prime size supported was $2^{32}-1$ but searches were generally done up to a max prime size of $\num{2000000}$. The PSLQ required 1024 bit precision numbers and the algorithm used was Bailey's algorithm m2 as specified in his paper \cite{Bailey2001} and ported fairly literally from his Fortran implementation.

The parts which needed to run fast were written in Go which is much easier to write than a lower level language like C or C++ but has most of the speed. It is also enables parallel code very simply. Python was used for all the data organization tasks, to extend the relations found with PSLQ and to double check the relations found.

The searches were run on a 32 core machine with 128GB of RAM over the period of a few months as the program evolved.

\section{Discussion}

\subsection{Lowest possible \texorpdfstring{$\lambda$}{Lehmer measure} for a relation with a given number of terms?}

Alferov's result gives us an insight into what the upper bound of the lowest Lehmer measure for a given number of terms might be.

\[
\lambda < \frac{3}{\log_{10} q_0}
\]

Alferov also proves that the number of terms $N$ in the sequence is

\[
N = O(q_0 \ln q_0)
\]

So

\[
N < O(q_0)
\]

Therefore the upper bound on the lowest Lehmer measure is

\[
\lambda < O(\frac{1}{\log_{10} N})
\]

Experiments indicate that using a constant of $1.03$  would cover the best relations found so far assuming that the best results for 7 terms and up are yet to be found.

\[
\lambda < \frac{1.03}{\log_{10} N}
\]

\subsection{Conclusion}

Finding the arctan relation with the lowest Lehmer measure for a given number of terms is an interesting mathematical and computational challenge. It seems likely that there are only a finite number of arctan relations of a given length. St\"ormer proved this for length 2 in \cite{Stormer1899} and Tomohiro Yamada proved this for length 3 in \cite{Yamada2025}.

We've shown new relations for 5 and 6 terms, and the 6 term relation has the lowest Lehmer measure of any relation known (other than the iterated relation \eref{alferov} with millions of digits).

The PSLQ based algorithm was very effective at finding short Machin-like formulae of up to 7 terms or so finding the existing and new best formulae for 5 terms and below. Combining PSLQ that with the idea of looking for $\acott{2}{x}$ terms produced the 6 and 7 term results with the lowest Lehmer measure. Attempting to extend the sequences with the floor-function iteration was very effective at finding formulae 8 or more terms and low Lehmer measure. It seems unlikely any of the 8 or more term formulae is the best possible but, with the exception of 9 terms where Wetherfield's 2009 relation (measure 1.4140) remains the best known, they are a big improvement on the previous best known.

\subsection{Future work}

There are several further directions for this work.

\begin{question}
  Explore more thoroughly arctan relations $\acott{a}{b}$ with $a\ne1$. Once found with PSLQ these relations could be expanded with the floor-function iteration rather than Todd's process as used for the $\acott{2}{b}$ cases here.
\end{question}

\begin{question}
  Explore more the selection of the prime group as the input to the PSLQ process. The methods used worked but didn't cover the whole space so there are likely more relations lurking.
\end{question}

\begin{question}
  Explore the question of whether there are a finite number of arctan relations of a given length. This could help put a bound on the largest prime needed for a given relation length.
\end{question}

\section{Data availability statement}

The source code and instructions on how to use it can be found on Github \cite{Github} at \url{https://github.com/ncw/find_arctan_formulae}.

The full results of this study including all Machin-like relations found are openly available in Zenodo \cite{Zenodo} at \url{https://doi.org/10.5281/zenodo.16758216}.

\section*{Declaration of interest}
The author declares no competing interests.

\section*{Changes in version 2}
\begin{itemize}
\item The previously best known 5, 6 and 9 term relations (Section~\ref{sec:best}, Table~\ref{table:terms1}) were wrong in version 1. The corrected relations come from Wetherfield and Chien-Lih's Machination lists, as preserved in the machin-like.org database \cite{MachinLike}. Wetherfield's 7 term relation dates from 2004, not 2003. The new 5 and 6 term relations reported here remain the best known.
\item The iteration in Section~\ref{sec:iteration} was attributed to Alferov. It was introduced by Abrarov, Siddiqui, Jagpal and Quine \cite{Abrarov2022}; Alferov's contributions are the nearest-integer variant, the termination proof and the bound on the Lehmer measure. Thanks to Sanjar Abrarov for pointing this out.
\item References to equation \eref{twoa} pointed at the wrong number in version 1.
\end{itemize}

% references

\bibliographystyle{plain}
\bibliography{refs} % References are in the refs.bib file

\begin{thebibliography}{10}

\bibitem{Abrarov2022}
Sanjar~M. Abrarov, Rehan Siddiqui, Rajinder~K. Jagpal, and Brendan~M. Quine.
\newblock A new form of the {Machin}-like formula for $\pi$ by iteration with
  increasing integers.
\newblock {\em Journal of Integer Sequences}, 25, 2022.
\newblock Article 22.4.5. Available at
  \url{https://cs.uwaterloo.ca/journals/JIS/VOL25/Abrarov/abrarov5.html}, also
  arXiv:2108.07718.

\bibitem{Abrarov2023}
Sanjar~M. Abrarov, Rehan Siddiqui, Rajinder~K. Jagpal, and Brendan~M. Quine.
\newblock A generalized series expansion of the arctangent function based on
  the enhanced midpoint integration.
\newblock {\em AppliedMath}, 3(2):395--405, 2023.
\newblock \url{https://doi.org/10.3390/appliedmath3020020}.

\bibitem{Alferov2023}
Oleg~S. Alferov.
\newblock A rapidly converging {Machin}-like formula for $\pi$, 2023.

\bibitem{Arndt2001}
J{\"o}rg Arndt and Christoph Haenel.
\newblock {\em Pi --- Unleashed}, chapter 16. $\pi$ Formula Collection, pages
  223--238.
\newblock Springer Berlin Heidelberg, Berlin, Heidelberg, 2001.

\bibitem{Bailey2001}
David Bailey and David Broadhurst.
\newblock Parallel integer relation detection: Techniques and applications.
\newblock {\em Mathematics of Computation}, 70(236):1719--1736, 2001.

\bibitem{Bailey2003}
David~H Bailey.
\newblock Some background on {Kanada}'s recent pi calculation.
\newblock {\em manuscript, 16-May}, 2003.

\bibitem{Chien1997}
Hwang Chien-Lih.
\newblock 81.22 more {Machin}-type identities.
\newblock {\em The Mathematical Gazette}, 81(490):120--121, 1997.

\bibitem{Zenodo}
Nick Craig-Wood.
\newblock {Machin-like} formulae for $\pi$ found using {PSLQ}, August 2025.
\newblock Available at \url{https://doi.org/10.5281/zenodo.16758216}.

\bibitem{Github}
Nick Craig-Wood.
\newblock Programs to calculate {Machin-like} formulae for $\pi$ using {PSLQ},
  2025.
\newblock Available at \url{https://github.com/ncw/find_arctan_formulae}.

\bibitem{PSLQ}
Helaman R.~P. Ferguson and David~H. Bailey.
\newblock A polynomial time, numerically stable integer relation algorithm.
\newblock Technical Report RNR-91-032, RNR, 1992.

\bibitem{Ferguson1999PSLQ}
Helaman R.~P. Ferguson, David~H. Bailey, and Steve Arno.
\newblock Analysis of {PSLQ}, an integer relation finding algorithm.
\newblock {\em Mathematics of Computation}, 68(225):351--369, 1999.

\bibitem{Lehmer1938}
Derrick~Henry Lehmer.
\newblock On arccotangent relations for $\pi$.
\newblock {\em The American Mathematical Monthly}, 45(10):657--664, 1938.

\bibitem{MachinLike}
Matthew~W. Scroggs.
\newblock machin-like.org: an encyclopedia of {Machin}-like formulae, 2025.
\newblock Available at \url{https://machin-like.org/}. Source data at
  \url{https://github.com/mscroggs/machin}, archived at
  \url{https://doi.org/10.5281/zenodo.17954248}.

\bibitem{Stormer1899}
Carl Störmer.
\newblock Solution compl\`ete en nombres entiers de l\textquoteright\'equation
  $m\arctan\frac{1}{x} + n\arctan\frac{1}{y} = k\frac{\pi}{4}$.
\newblock {\em Bulletin de la Soci\'et\'e Math\'ematique de France},
  27:160--170, 1899.

\bibitem{HackersDelight}
Henry~S. Warren\, Jr.
\newblock {\em Hacker's Delight}.
\newblock Addison-Wesley Professional, 2nd edition, 2012.

\bibitem{Wetherfield1996}
Michael Wetherfield.
\newblock The enhancement of {Machin}'s formula by {Todd}'s process.
\newblock {\em The Mathematical Gazette}, 80(488):333--344, 1996.

\bibitem{Wetherfield1997}
Michael Wetherfield.
\newblock 81.23 {Machin} revisited.
\newblock {\em The Mathematical Gazette}, 81(490):121--123, 1997.

\bibitem{Machination}
Michael Wetherfield and Hwang Chien-lih.
\newblock Computing pi: Lists of {Machin}-type (inverse cotangent) identities
  for pi/4, 2013.
\newblock The original website no longer exists; a partial copy is archived at
  \url{https://web.archive.org/web/20240204042153/http://www.machination.eclipse.co.uk/}
  and its relations are included in \cite{MachinLike}.

\bibitem{Yamada2025}
Tomohiro Yamada.
\newblock Three-term {Machin}-type formulae, 2025.

\end{thebibliography}

% Tables (each on its own page)

\renewcommand{\arraystretch}{1.5}

\clearpage
\begin{table}[ht]
  \centering\scriptsize
    \caption{The best known relations: 1-11 terms (found in this work unless marked)} \label{table:terms1}
      \begin{tabularx}{\linewidth}{c S[round-mode=places,round-precision=4,group-minimum-digits=5] Y}
      \toprule
      Terms & $\lambda$ & Relation for $\piF$ ($[x]=\acotx{x}$) \WHO{If discovered previously}\\
      \midrule
2 & 1.8511276523 & 4\ACOT{5}-1\ACOT{239} \WHO{Machin 1706}\\
3 & 1.7866075340 & 12\ACOT{18}+8\ACOT{57}-5\ACOT{239} \WHO{Gauss 1863}\\
4 & 1.5860413586 & 44\ACOT{57}+7\ACOT{239}-12\ACOT{682}+24\ACOT{12943} \WHO{St\"ormer 1896}\\
5 & 1.4571904502 & 29\ACOT{68}+42\ACOT{117}-15\ACOT{2675143}-5\ACOT{2976163}-5\ACOT{302342643} \\
6 & 1.3291269825 & 83\ACOT{107}+17\ACOT{1710}-44\ACOT{225443}-68\ACOT{2513489}+22\ACOT{42483057}+34\ACOT{7939642926390344818} \\
7 & 1.3408464538 & 83\ACOT{107}+17\ACOT{1710}-22\ACOT{103697}-24\ACOT{2513489}-44\ACOT{18280007883}+12\ACOT{7939642926390344818}+22\ACOT{3054211727257704725384731479018} \WHO{Wetherfield 2004}\\
8 & 1.4167155254 & 83\ACOT{107}+17\ACOT{1799}+29\ACOT{110443}+22\ACOT{4841182}+17\ACOT{5675477}-17\ACOT{230921399215798}-17\ACOT{180288246462881610746480172447}+17\ACOT{68258088806587561619922885741255927889500440723402481904848} \\
9 & 1.4140461772 & 83\ACOT{107}+17\ACOT{1710}-44\ACOT{226043}-68\ACOT{2513489}-22\ACOT{109027476193}+22\ACOT{5774883029189818}+34\ACOT{7939642926390344818}+44\ACOT{3375905320682366575989}-22\ACOT{19237152089720411315324480876866814545312307000809567958779407318} \WHO{Wetherfield 2009}\\
10 & 1.4419042849 & 83\ACOT{107}+17\ACOT{1710}-22\ACOT{112917}-68\ACOT{2513489}+22\ACOT{122261817}-22\ACOT{105765955954736771}+34\ACOT{7939642926390344818}+22\ACOT{37315951755534592843736520276224678}-22\ACOT{5183120950735156251347952259333584038201011704144309250777626572093950}-22\ACOT{483565370219094780653589520268389892857600447041019458360871023697086072184994990755122112972383385497181774829503910142319319126350706751068} \\
11 & 1.4406475833 & 83\ACOT{107}+17\ACOT{1710}-44\ACOT{225768}-68\ACOT{2513489}+22\ACOT{92866807}-22\ACOT{25006102376411521}+34\ACOT{7939642926390344818}+22\ACOT{1456552709819752015648621402178712}-22\ACOT{13519446218663109947650509037037237559876973237968142522118037123610}-22\ACOT{991955953079245266974385548959755309175792978557592484113060221474505251076876235206434243729747213526161424285017818981041750239211773}+22\ACOT{10134959112348344527804392569520494971086697287357393590158406211657079084503520349173141299677656646171095816120989447418951150541424947997810136848062107605261111776545978270751363275882456888729201312941418105672728754987083678664924509436369478551407903599805359788132} \\
      \bottomrule
    \end{tabularx}
\end{table}

\clearpage
\begin{table}[ht]
  \centering\scriptsize
    \caption{The best relations found: 12-15 terms} \label{table:terms2}
      \begin{tabularx}{\linewidth}{c S[round-mode=places,round-precision=4,group-minimum-digits=5] X}
      \toprule
      Terms & $\lambda$ & Relation for $\piF$ ($[x]=\acotx{x}$) \\
      \midrule
12 & 1.4921095533 & 83\ACOT{107}+17\ACOT{1710}-88\ACOT{452761}-68\ACOT{2513489}-22\ACOT{75757291}+22\ACOT{26619138353848472}+44\ACOT{46406309881072682}+34\ACOT{7939642926390344818}+22\ACOT{1552740245358442879545096003111139}+22\ACOT{7892397331196677585383379835983113026631455870097764779488866226405}-22\ACOT{65699296015786320869593856066716893749129800185319202819151454852291380141816241738912306404253371756817487701294327908830743814165355}+22\ACOT{8632794993939832670748112451230064842497999893667491222538502729787729319898264659861783713136001487066117764056406191664287179911376869841893791277365861124913985478061814254780518073614587140507633566355499913511438490258926207106663273692310676785272596254378717407} \\
13 & 1.6259278468 & 83\ACOT{107}+17\ACOT{1710}-22\ACOT{103615}-24\ACOT{2513489}-22\ACOT{64346526}+44\ACOT{86718193}-22\ACOT{8756857644969543}+12\ACOT{7939642926390344818}-22\ACOT{326062356717297569844818}+22\ACOT{1867355102894027035134301488735408}+22\ACOT{19826171456587659909082231353393318806646988224831087078872555617922}+22\ACOT{917179840794026255056120358332428626619503750708425137609940317818147293581189653664388847143573874363412854773508297216719361688777610}+22\ACOT{2523656581076866043300463048811818630190963714279867338534700369442654018006950072819452024396983318050290174370560841720378277657606609331173829131728776326166279365087406454460617465557398462538631814699502188303896290489771688980525226763109155846348268367686423158693} \\
14 & 1.6764894343 & 1345\ACOT{1710}+166\ACOT{11654}-581\ACOT{48443}-708\ACOT{225443}-732\ACOT{2513489}+354\ACOT{42483057}+1162\ACOT{54252061}-83\ACOT{1558810172}+83\ACOT{5147177440427785908}+366\ACOT{7939642926390344818}-581\ACOT{79839669135881505029582}+83\ACOT{222969562996807089602885876643155659684}-83\ACOT{58938649855306003782663801408220460798195205953191133428534614497951030866573}+83\ACOT{11810799119005632230130264080871074758919453586112153180108213002548927458013086740177685105609587262830477714762294167394949437703240807485402361548665295} \\
15 & 1.8061930072 & 83\ACOT{239}+764\ACOT{1710}-332\ACOT{48443}-752\ACOT{452761}-400\ACOT{2513489}+664\ACOT{54252061}-188\ACOT{75757291}+188\ACOT{26619138353848472}+376\ACOT{46406309881072682}+200\ACOT{7939642926390344818}-332\ACOT{79839669135881505029582}+188\ACOT{1552740245358442879545096003111139}+188\ACOT{7892397331196677585383379835983113026631455870097764779488866226405}-188\ACOT{65699296015786320869593856066716893749129800185319202819151454852291380141816241738912306404253371756817487701294327908830743814165355}+188\ACOT{8632794993939832670748112451230064842497999893667491222538502729787729319898264659861783713136001487066117764056406191664287179911376869841893791277365861124913985478061814254780518073614587140507633566355499913511438490258926207106663273692310676785272596254378717407} \\
      \bottomrule
    \end{tabularx}
\end{table}

\end{document}